\documentclass[11pt]{amsart}
\title[Schur tensor product]{
Schur tensor product of operator spaces
}
\usepackage[ansinew]{inputenc}
\usepackage{yfonts}
\usepackage{calligra}
\usepackage{amscd}
\usepackage{amssymb, latexsym}
\usepackage{mathrsfs}
\usepackage{fancyhdr}
\usepackage{enumerate}

\usepackage{textcomp, phonetic}
\usepackage{setspace}
\usepackage{latexsym}
\newcommand{\ciii}{\stackrel{cb}{=}}
\usepackage[all]{xy}
\usepackage{tikz}
\theoremstyle{plain}
\newtheorem{thm}{\sc Theorem}[section]
\newtheorem{cor}[thm]{\sc Corollary}
\newtheorem{lem}[thm]{\sc Lemma}

\newtheorem{prop}[thm]{\sc Proposition}

\newtheorem{defi}[thm]{\sc Definition}

\newcommand{\ot}{\otimes}

\newcommand{\os}{\otimes^s}
\newcommand{\C}{\mathbb{C}}

\newcommand{\p}{\varphi}
\newcommand{\pb}{\overline{\varphi}}
\newcommand{\ds}{\displaystyle}
\newcommand{\al}{\alpha}
\newcommand{\bt}{\beta}
\newenvironment{pf}{\noindent {\sc Proof:}}{\hfill $\Box$}
\begin{document}

\author[V. Rajpal, A. Kumar And T. Itoh]{
VANDANA RAJPAL, AJAY KUMAR$^{1}$ AND TAKASHI ITOH}
\address{Department of Mathematics\\
University of Delhi\\
Delhi\\
India.}
\email{vandanarajpal.math@gmail.com}
\address{Department of Mathematics\\
University of Delhi\\
Delhi\\
India.}
\email{akumar@maths.du.ac.in}
\address{Department of Mathematics\\
University of Gunma\\
Gunma\\
Japan.}
\email{itoh@gunma-u.ac.jp}
\thanks{}
\footnotetext[1]{- corresponding author}
\maketitle
\begin{abstract}
We develop a systematic study  of the schur tensor product both in the category of operator spaces and in that of
$C^*$-algebras.
\end{abstract}

\section{Introduction}
An operator space is a closed subspace  of the  space $B(H)$ of all bounded operators on a Hilbert space $H$. The fundamental and systematic developments in the theory of  tensor product of operator spaces have been evolved considerably,  see e.g.  ~\cite{blecher}, ~\cite{ER}. In this category, the Haagerup tensor norm is the natural one for compatibility with  the continuity of the completely bounded bilinear maps, and the operator space projective tensor norm is for the jointly completely  bounded  bilinear maps.
For operator spaces $V$ and $W$,  and elements  $x=[x_{ij}]\in M_n(V)$ and $y=[y_{ij}]\in M_n(W)$, we define an element
 $x \circ y \in M_n(V \otimes W) $ by $x \circ y=[x_{ij}\otimes y_{ij}]$. Note that $x \circ y=[e_{11},\; e_{22},\; e_{33},\cdots, e_{nn}] (x\ot y)[e_{11},\; e_{22},\; e_{33},\cdots, e_{nn}]^{t} $, where $x\ot y$ denotes the kronecker tensor product, and $\{e_{ij}\}$ are the standard basis of $M_n$, $n\in \mathbb{N}$. Each element $u$ in $M_p(V\ot W)$, $p\in \mathbb{N}$, can be written as $u=\alpha(x \circ y)\beta$ for some $x\in M_n(V)$,  $y\in M_n(W)$, $\alpha\in M_{p,n}$, and $\beta\in M_{n,p}$, $n\in \mathbb{N}$, and we define
\[\|u\|_s=\inf\{\|\alpha\|\|x\|\|y\|\|\beta\|\}\]
where infimum is taken over arbitrary decompositions as above.  Let $V \otimes_s W=(V \otimes W, \|\cdot\|_s)$, and define  the schur tensor product $V \otimes^s W$ to be the completion of $V \otimes W$ in this norm. In Section 2 we look at the schur  tensor norm in the
context of operator spaces, and  show that this is an operator space matrix norm.
 We  also introduce   the notion of schur bounded bilinear maps and show that the schur tensor norm may be used
 to linearize them. This parallel development of the theory of schur tensor product of operator spaces will play a vital role in the theory of
 operator spaces. In analogy to the operator space projective tensor product, the schur tensor product turns out to be commutative, projective,
  and functorial however we don't know whether is it associative or not. At the end of this section, we define a new tensor norm,
  which we denote by $\|\cdot\|_{s'}$, and it will be seen  that the dual of the schur tensor norm  is in fact  $\|\cdot\|_{s'}$-norm.
 Section 3 is concerned with the equivalence of the schur tensor norm with the various other norms.

Recall that the operator space projective tensor norm on the algebraic tensor product of two operator spaces $V$ and $W$ is defined as, for $u\in V\otimes W$, $\|u\|_{\wedge}=\inf\{\|\alpha\|\|x\|\|y\|\|\beta\|\},$ the infimum is taken over $p,q\in \mathbb{N}$ and all the ways to write $u=\alpha (x\ot y)\beta$, where $\alpha\in M_{1,pq}$, $\beta\in M_{pq,1}$, $x\in M_{p}(V) $ and $y\in M_{q}(W)$, and
$x\otimes y=(x_{ij}\otimes y_{kl})_{(i,k),(j,l)}\in M_{pq}(V\otimes W)$. The operator space projective  tensor product $V\widehat{\otimes} W$ is defined to be the completion of $V\otimes W$ in the norm $\|\cdot\|_{\wedge}$~\cite{blecher}. The jointly completely bounded  norm of a bilinear map $\phi: V\times W\to Z$ is defined to be $\|\phi\|_{jcb}=\sup\{\|\phi_{(n)}\|:n\in\mathbb{N}\}$, where $\phi_{(n)}: M_n(V)\times M_n(W)\to M_{n^2}(Z)$ is  given by $\phi_{(n)}([v_{ij}],[w_{kl}])=(\phi(v_{ij},w_{kl}))$.

\section{Schur tensor product of operator spaces}
\begin{thm}
For operator spaces $V$ and $W$, $\|\cdot \|_s$ is an operator space matrix norm on $V \otimes W$.
\end{thm}
\begin{pf}
%
Given $u_1\in V \otimes W$,  $u_2\in V \otimes W$ and $\epsilon > 0$, choose $\alpha_1\in M_{1,r}$, $x_1\in M_r(V)$, $y_1\in M_r(W)$, $\beta_1\in M_{r,1}$, and $\alpha_2\in M_{1,p}$,  $x_2\in M_p(V)$, $y_2\in M_p(W)$, $\beta_2\in M_{p,1}$, $p,r\in \mathbb{N}$ such that $u_1=\alpha_1(x_1\circ y_1)\beta_1$, $u_2=\alpha_2(x_2\circ y_2)\beta_2$ and  $\|\alpha_1\| \|x_1\|\|y_1\|\|\beta_1\|< \|u_1\|_{s}+\epsilon$, $\|\alpha_2\| \|x_2\|\|y_2\|\|\beta_2\|< \|u_2\|_{s}+\epsilon$, we may assume  that $\|x_i\|=\|y_i\|=1$ and $\|\alpha_i\|=\|\beta_i\|\leq (\|u_i\|_s+\epsilon)^{\frac{1}{2}}$, for $i=1,2$. Let $\alpha=[\begin{array}{cc}
                                                                                   \alpha_1 & \alpha_2
                                                                                 \end{array}]
$, $\beta=[\begin{array}{cc}
                                                                                   \beta_1 & \beta_2
                                                                                 \end{array}]^{t}
$, $v:= x_1\oplus x_2$, and $w:= y_1\oplus y_2$. Then
$v \circ w=\left[\begin{smallmatrix}
               x_1 & 0 \\
               0 & x_2  \\
  \end{smallmatrix}\right]\circ \left[\begin{smallmatrix}
               y_1 & 0 \\
               0 & y_2  \\
  \end{smallmatrix}\right]= \left[\begin{smallmatrix}
               x_1 \circ y_1 & 0 \\
               0 & x_2 \circ y_2  \\
  \end{smallmatrix}\right]$,
and so $u_1 + u_2=\alpha  (v \circ w)\beta$.
Now, by Ruan's axioms of operator spaces~\cite{effros}  and the $C^*$-identity, we have \\
  \hspace*{ 4 cm} $\|u_1+ u_2\|_s\leq  \|\alpha\|\|\beta\|$,\\
  \hspace*{ 5.8 cm} $\leq \frac{1}{2} (\|\alpha\|^2+\|\beta\|^2)$,\\
   \hspace*{ 5.8 cm} $=\frac{1}{2} \|\left[\begin{smallmatrix}
   \al_1  & \al_2 \\
               0 & 0  \\
  \end{smallmatrix}\right]\|^2 +\|\left[\begin{smallmatrix}
   \beta_1  & 0 \\
               \beta_2 & 0  \\
  \end{smallmatrix}\right]\|^2$\\
   \hspace*{ 5.8 cm} $=\frac{1}{2} \|\left[\begin{smallmatrix}
   \al_1  & \al_2 \\
               0 & 0  \\
  \end{smallmatrix}\right]\left[\begin{smallmatrix}
   \al_{1}^*  & 0 \\
               \al_{2}^* & 0  \\
  \end{smallmatrix}\right]\|+\|\left[\begin{smallmatrix}
   \bt_{1}^*  & \bt_{2}^* \\
               0 & 0  \\
  \end{smallmatrix}\right]\left[\begin{smallmatrix}
   \bt_{1}  & 0 \\
               \bt_{2} & 0  \\
  \end{smallmatrix}\right]\|$\\
 \hspace*{ 5.8 cm} $= \frac{1}{2}\|\left[\begin{smallmatrix}
   \al_{1}\al_{1}^*+\al_{2}\al_{2}^*  & 0 \\
               0 & 0  \\
  \end{smallmatrix}\right]\|+\|\left[\begin{smallmatrix}
   \bt_{1}^*\bt_{1}+\bt_{2}^*\bt_{2}  & 0 \\
               0 & 0  \\
  \end{smallmatrix}\right]\|$\\
 \hspace*{ 5.8 cm} $= \frac{1}{2} \|\alpha_1 \al_{1}^* +\al_2\al_2 ^{*}\|+\| \bt_{1}^* \bt_1+\bt_2 ^{*}\bt_2\|$,\\
  \hspace*{ 5.8 cm} $\leq \frac{1}{2} (\|\alpha_1\|^2 +\|\al_2 \|^2+\|\beta_1\|^2+\|\beta_2\|^2)$,\\
  \hspace*{ 5.8 cm} $\leq \|u_1\|_s +\|u_2 \|_s+ 2 \epsilon$. \\
Since $\epsilon>0$ is arbitrary, so we have the subadditivity. For any scalar $c\in \mathbb{C}$, it is clear that $\|cu\|_s=|c|\|u\|_s$.

Let $u=\al (x \circ y)\beta\in V\ot W$, $\al=[\al_i]\in M_{1,p}$, $\beta=[\beta_j] \in M_{p,1}$, $x\in M_p(V)$, and $y\in M_p(W)$, we can write $u$ as  $u=\left(\begin{matrix}
               \al_1 \; 0 \;\ldots 0 \; 0\; \al_2\ldots \; 0\ldots 0\;0\;\ldots \al_p
  \end{matrix}\right)(x\ot y)$
  $ \left(\begin{matrix}
               \beta_1  0 \ldots 0 \;\; 0 \;\beta_2\ldots  0\;\ldots \;0\; 0 \ldots \beta_p
  \end{matrix}\right)^{t}$, which is a representation in the set $\{\gamma(x\ot y)\lambda: \gamma\in M_{1,p^{2}}, \lambda\in M_{p^{2},1}, x\in M_p(V), y\in M_p(W) \}$.  Again, by Ruan's axioms and the $C^*$-identity,
$\|\left(\begin{matrix}
               \al_1 \; 0 \;\ldots \;0 \; 0\; \al_2\ldots \; 0\ldots \;0\;0\;\ldots\; \al_p
  \end{matrix}\right)\| =\|\al\|$ and $\|\left(\begin{matrix}
               \beta_1 \; 0\; \ldots \;0 \; 0\; \beta_2\;\ldots \; 0\;\ldots \;0\;0\;\ldots \;\beta_p
  \end{matrix}\right)^{t}\|=\|\beta\|$. Hence $\|u\|_\wedge \leq  \|u\|_s$, which shows that  $\|\cdot\|_s$ is a norm on $V \otimes W $.

Finally, to see the Ruan's axioms. Suppose that $u_1\in M_m(V \otimes W)$,  $u_2\in M_n(V \otimes W)$ and $\epsilon > 0$, we may find decompositions, for $i=1,2$, $u_i=\alpha_i(x_i \circ y_i)\beta_i$
with $\|x_i\|=\|y_i\|=1$ and $\|\alpha_i\|=\|\beta_i\|\leq (\|u_i\|_s+\epsilon)^{\frac{1}{2}}$ as above, where $\alpha_1\in M_{m,r}$, $\alpha_2\in M_{n,s}$,  $\beta_1\in M_{r,m}$, and $\beta_2\in M_{s,n}$.
Let $v:= x_1\oplus x_2$ and $w:= y_1\oplus y_2$. Then $u_1 \oplus u_2=\left[\begin{smallmatrix}
               \alpha_1  & 0 \\
               0 & \alpha_2   \\
  \end{smallmatrix}\right] (v \circ w)\left[\begin{smallmatrix}
               \beta_1  & 0 \\
               0 & \beta_2   \\
  \end{smallmatrix}\right]$. So, by Ruan's axioms for $V$ and $W$, we have  \\
  \hspace*{ 3 cm} $\|u_1 \oplus u_2\|_s\leq  \big\|\left[\begin{smallmatrix}
               \alpha_1  & 0 \\
               0 & \alpha_2   \\
  \end{smallmatrix}\right]\big\| \big\|  \left[\begin{smallmatrix}
               \beta_1  & 0 \\
               0 & \beta_2   \\
  \end{smallmatrix}\right]\big\|$.\\
Let $t=\max\{m+n, r+s\}$. By adding rows and columns of zeros to matrices $\left[\begin{smallmatrix}
               \alpha_1  & 0 \\
               0 & \alpha_2   \\
  \end{smallmatrix}\right]$ and $\left[\begin{smallmatrix}
               \beta_1  & 0 \\
               0 & \beta_2   \\
  \end{smallmatrix}\right]$, we obtain  the new matrices, say, $T_1\in M_{t}(\mathbb{C})$ and $T_2\in M_{t}(\mathbb{C})$ such that $\left[\begin{smallmatrix}
               \alpha_1  & 0 \\
               0 & \alpha_2   \\
  \end{smallmatrix}\right]$ and $\left[\begin{smallmatrix}
               \beta_1  & 0 \\
               0 & \beta_2   \\
  \end{smallmatrix}\right]$  are the upper left hand corner of $T_1$ and $T_2$, respectively, and so applying the $C^*$-identity,  we have\\
  \hspace*{ 2 cm} $\|u_1 \oplus u_2\|_s\leq  \big\|\left[\begin{smallmatrix}
               \alpha_1 \alpha_{1}^{*} & 0 \\
               0 & \alpha_2 \alpha_{2}^{*}  \\
  \end{smallmatrix}\right]\big\|^{1/2}\big\|\left[\begin{smallmatrix}
               \beta_1\beta_{1}^{*}  & 0 \\
               0 & \beta_2\beta_{2}^{*}   \\
  \end{smallmatrix}\right]\big\|^{1/2}$\\\\
  \hspace*{3.8 cm} $= \|\alpha_1 \alpha_{1}^{*} \oplus \alpha_2 \alpha_{2}^{*}\|^{1/2} \|\beta_1\beta_{1}^{*} \oplus \beta_2\beta_{2}^{*}\|^{1/2}$\\\\
  \hspace*{3.8 cm} $=(\max\{\|\alpha_1\|^{2},\|\alpha_2\|^{2}  \})^{1/2}(\max\{\|\beta_1\|^{2},\|\beta_2\|^{2}  \})^{1/2}$ \\
\hspace*{3.8 cm} $\leq \max\{\|u_1\|_s, \|u_2\|_s\}+\epsilon$.\\
Since $\epsilon$ is arbitrary, so  $\|u_1 \oplus u_2\|_s\leq \max\{\|u_1\|_s, \|u_2\|_s\}$. Now let $\gamma \in M_{p, m}$ and $\delta \in M_{m, p}$, then $\gamma u_1 \delta= (\gamma \alpha_1)(x_1 \circ y_1)(\beta_1 \delta)$, and so
$\|\gamma u_1 \delta\|_s\leq \|\gamma \alpha_1\| \|\beta_1 \delta\| \leq \|\gamma \| \| \delta\| (\|u_1\|_s +\epsilon )$. Because $\epsilon$ was arbitrary,  $\|\gamma u_1 \delta\|_s\leq   \|\gamma \| \| \delta\| \|u_1\|_s.$  Hence the result follows from (~\cite{effros}, Proposition 2.3.6).
\end{pf}

\begin{thm}\label{sc13}
For operator algebras $X$ and $Y$, $X\otimes^{s}Y$ is a Banach algebra, and it is $^*$-algebra provided both $X$ and $Y$ have isometric involution. Furthermore, if $X$ and $Y$ are approximately unital then $X\otimes^{s}Y$ has  a bounded approximate identity.
\end{thm}
\begin{pf}
Let $u$, $v\in X\otimes Y$ with $u=\sum_{i,j=1}^{n}\alpha_i (x_{ij} \otimes y_{ij})\beta_j$, and $v=\sum_{k,l=1}^{m}\gamma_k (z_{kl} \otimes w_{kl})\delta_l$. Then $uv=\sum_{i,j,k,l}\alpha_i \gamma_k(x_{ij}z_{kl} \otimes y_{ij}w_{kl})\beta_j\delta_l,$
which can be further written as $uv=[\alpha_i \gamma_k][x_{ij}z_{kl}] \circ [y_{ij}w_{kl}][\beta_j\delta_l]$,
and so\\
\hspace*{3 cm}$\|uv\|_s\leq \|[\alpha_i \gamma_k]\|\|[x_{ij}z_{kl}]\| \| [y_{ij}w_{kl}]\|\|[\beta_j\delta_l]\|$.\\
Thus,  by (~\cite{effros}, Proposition 2.1.1) and   the fact that the operator algebras are completely contractive, we have \\
\hspace*{3 cm} $\|uv\|_s\leq\|[\alpha_i]\|\|[ \gamma_k]\|\|[x_{ij}\|\|[z_{kl}]\| \| [y_{ij}]\|\|[w_{kl}]\|\|[\beta_j]\|\|[\delta_l]\|.$\\
Hence, $\|uv\|_s\leq \|u\|_s\|v\|_s$ for all $u,v\in X \otimes_s Y$. So we may extend the product on $X\otimes_s Y$ to $X \otimes^s Y$, so that
$X\otimes^{s}Y$ is a Banach algebra.

For the $*$-part, let  $u\in X\otimes Y$ with $u=\alpha (x \circ y)\beta$ then $u^*=\beta^* (x^* \circ y^*)\alpha^*$. So $\|u^*\|_s\leq \|\beta^* \| \|x^*\|\| y^*\|\|\alpha^*\|= \|\beta\| \|x\|\| y\|\|\alpha\|$. Thus $\|u^*\|_s\leq \|u\|_s$. Similarly, $\|u\|_s\leq \|u^*\|_s$.

One can easily verify that $\|\cdot\|_s \leq \|\cdot\|_\gamma$ on $X \otimes Y $, giving that  $\|\cdot\|_s$ is an admissible cross norm on $X \otimes Y$. Therefore, $X\otimes^{s}Y$ has a bounded approximate identity, whenever $X$ and $Y$ are approximately unital.
\end{pf}

In particular, for $C^{*}$-algebras $A$ and $B$, $A\otimes^{s}B$ is a Banach $^*$-algebra with a bounded approximate identity, and it is a  $C^*$-algebra if and only if either $A=\mathbb{C}$ or $B=\mathbb{C}$, which follows directly by (~\cite{blecher}, Corollary 5.13) and (~\cite{geom}, Corollary 3). Also, (~\cite{laur}, Theorem 4.8) yields that the enveloping $C^*$-algebra, $C^*(A\otimes^{s}B)$,  of  $A\otimes^{s}B$ is $A\otimes_{\max}B$, the maximal tensor product of $A$ and $B$.

As with the Haagerup tensor norm and the operator space projective tensor norm, we define an intimately related class of bilinear maps for the schur tensor product of operator spaces  by mimicking the schur multiplication of matrices.
\begin{defi}
Given operator spaces $V$, $W$ and $Z$, a bilinear map $\varphi: V \times W \to Z$ is said to be schur bounded bilinear map if the associated maps $\varphi_{n}:M_n(V)\times M_{n}(W)\to M_{n}(Z)$ given by
\begin{equation}
\varphi_n \big( (v_{ij}), (w_{ij}) \big )=\big(\varphi(v_{ij},w_{ij})\big), \, \, n \in \mathbb{N}\notag
\end{equation}
are uniformly bounded, and in this case we denote $\|\varphi\|_{sb} = \displaystyle\sup_{n}\|\varphi_n\|$.
\end{defi}
Equivalently, a bilinear map $\varphi: V \times W \to Z$ is schur bounded if and only if the associated maps $\varphi_n:M_n(V)\times M_{n}(W)\to M_{n}(Z)$ given by
\begin{equation}
\varphi_n (\al \ot v, \bt \ot w)= \al \circ \bt \ot \varphi(v,w), \al, \bt\in M_n, v\in V, w\in W,\notag
\end{equation}
where $\al \circ \bt$ is the usual schur multiplication of matrices, are uniformly bounded. Indeed,  for $\al=[\al_{ij}]$ and $\bt=[\bt_{ij}]$, $\al \ot v=[\al_{ij}v]$ and $\bt \ot w=[\bt_{ij}w]$ by using the identification $M_n(V) \cong M_n \ot V$. Thus, by the above definition, $\varphi_n (\al \ot v, \bt \ot w)=\varphi_n ([\al_{ij}v],[\bt_{ij}w] )= \big(\varphi(\al_{ij}v,\bt_{ij}w)\big)= (\al_{ij}\bt_{ij}\varphi(v,w))= [\al_{ij}\bt_{ij}]\ot \varphi(v,w)= \al \circ \bt \ot \varphi(v,w)$. Also, any jointly completely bounded bilinear map  is schur bounded. This is immediate from the relation $\phi_{n}([x_{ij}],[y_{ij}])= (\phi(x_{ij},y_{ij}))=\alpha \phi_{(n)}([x_{ij}],[y_{kl}])\beta$ for $\alpha=[e_{11}, e_{22},\dots, e_{nn}]$ and $\beta=\alpha^{t}$.

It is easy to see that  $\|\cdot \|_{sb}$ is a norm on the linear space $SB(V \times W , Z)$, the space of all schur bounded bilinear maps.
We now show that the schur tensor norm linearizes the schur bounded bilinear maps, that is $ (V\os W)^*= SB(V\times W, \C). $

\begin{prop}\label{sc8}
If $V$, $W$ and $X$ are operator spaces, then there is a natural  isometric identification
\[CB(V\os W, X)= SB(V\times W, X).\]
\end{prop}
\begin{pf}
Let $\varphi: V \times W \to X$ be a schur bounded  bilinear map. Then there exists a unique linear mapping $\pb:  V \otimes W\to X$ such that $\pb(x\ot y)= \p(x,y)$ for all $x\in V$ and $y\in W$. For $u=\alpha(v\circ w)\beta\in M_n(V \otimes W)$, $\al\in M_{n,p}$, $v\in M_p(V)$, $w\in M_p(W)$, and $\beta\in M_{p,n}$,  we have $\pb_n(u)=  \alpha \p_{p}([v_{ij}],[w_{ij}])\beta$,\\
and so\\
\hspace*{3 cm} $\| \pb_n(u)\|\leq \|\alpha\|\|\p_{p}\|\|[v_{ij}]\|\|[w_{ij}]\|\|\beta\|$,\\
\hspace*{4.4259 cm} $\leq \|\p\|_{sb}  \|\alpha\|\|[v_{ij}]\|\|[w_{ij}]\|\|\beta\|$.\\
Since the above is true for every representation of $u$, so it follows that  $\|\pb\|_{cb}\leq  \|\p\|_{sb} $. For the converse part, note that  $\p_{p} ([v_{ij}],[w_{ij}])= (\p(v_{ij},w_{ij}))=\pb_{p}(v \circ w)$, for $v=[v_{ij}]$ and $w=[w_{ij}]$. So  $\|\p_{p} ([v_{ij}],[w_{ij}])\|=\|\pb_{p}(v \circ w)\|\leq\|\pb\|_{cb}\|v \circ w\|_s\leq\|\pb\|_{cb}\|v\| \|w\|$. Hence it follows that $\|\p\|_{sb}\leq \|\pb\|_{cb}$. Also the operator $\pb : V \otimes_s W\to X $ has a unique extension to an operator $\pb : V \otimes^s W\to X $ with the same norm.

To see the surjectivity, let $\psi\in CB(V\os W, X)$. We can define a blinear map $\p: V \times W\to X$ such that  $\p(v,w)=\psi(v\ot w)$. Using the same calculations as done in the above paragraph just replacing $\pb$ by $\psi$, we obtain  the required result.
\end{pf}

The  above identification yields a new formula for the schur tensor norm:
\[\|u\|_s=\sup\{|\p(u)|: \p\in SB(V\times W, \C), \|\p\|_{sb}\leq 1\}.\]

The following corollary, which shows that the schur tensor product is functorial, is a simple consequence of the above
Proposition.
\begin{cor}\label{sc113}
Let  $V $, $V_1$, $W$ and $W_1$ be operator spaces. Given the complete contractions $\varphi: V\to V_1$ and $\phi: W\to W_1$, the corresponding map
$\varphi \otimes \phi: V\otimes W\to  V_1\otimes W_1$ extends to a complete contraction map $\varphi \otimes^s \phi: V\otimes^s W\to  V_1\otimes^s W_1$.
\end{cor}
%

For operator spaces $V$ and $W$, let $B(V,W)$ denote the bounded linear maps from $V$ to $W$. Define the matrix norm structure on $B(V,W)$  by identifying $[f_{ij}]\in M_n(B(V,W))$ with the map $[f_{ij}]:M_n(V)\to M_n(W)$   defined by $[f_{ij}]([x_{ij}])=[f_{ij}(x_{ij})]$. Let $_{s}B(V,W)$ denote  the space $B(V,W)$ with this matrix norm structure. Then we have the following:

\begin{prop}
For operator spaces $V$, $W$ and $Z$,  $SB(V \times W, Z)=CB(V,_{s}B(W,Z))$, here equal sign signify the isometric isomorphism.
\end{prop}
\begin{pf}
Let $u: V \times W \to Z$ be a schur  bounded bilinear map, write $\widetilde{u}$ for the map from $V$ to the set of functions from $W$ to $Z$ defined by
$\widetilde{u}(v)(w) = u(v, w)$, $v\in V$ and $w\in W$. Then $\|\widetilde{u}\|_{cb}=\sup\{\|[\widetilde{u}(v_{ij})]\|_{M_n(_{s}B(W,Z))}: [v_{ij}]\in M_n(V)_1, n\in \mathbb{N}\}= \sup\{\|[\widetilde{u}(v_{ij})(w_{ij})]\|: [v_{ij}]\in M_n(V)_1, [w_{ij}]\in M_n(W)_1, n\in \mathbb{N}\}=\sup\{\|[u(v_{ij},w_{ij})]\|: [v_{ij}]\in M_n(V)_1, [w_{ij}]\in M_n(W)_1, n\in \mathbb{N}\}=\|u\|_{sb}$. For the converse, let $v\in CB(V,_{s}B(W,Z))$ and define $u(x, y) = v(x)(y)$, then reversing the last argument shows that $u$ is schur  bounded.
\end{pf}


The proof of the following propositions are essentially the same  as those for the operator space projective tensor product, so we skip them.

\begin{prop}\label{sc11}\emph{[Projective]}
Let  $V $, $V_1$, $W$ and $W_1$ be operator spaces. Given the complete quotient maps $\varphi: V\to V_1$ and $\phi: W\to W_1$, the corresponding map
$\varphi \otimes \phi: V\otimes W\to  V_1\otimes W_1$ extends to a  complete quotient map $\varphi \otimes^s \phi: V\otimes^s W\to  V_1\otimes^s W_1$. \\
Furthermore, \[\ker \varphi \otimes^s \phi=cl\{\ker \varphi \otimes W+ V\ot \ker\phi\}.\]
\end{prop}

\begin{prop}\label{sc1}\emph{[Symmetric]}
Given operator spaces $V$, $W$ and $Z$, we have completely isometric isomorphism:
\[V \otimes^s W \ciii W\otimes^s V.\]
\end{prop}
%

The next proposition gives the general representation of an element of the schur tensor product.
\begin{prop}\label{sc4}
Given operator spaces $V$ and $W$, if $u\in M_n(V \otimes^s W)$ then
\[\|u\|_s=\inf\{\|\alpha\|\|x\|\|y\|\|\beta\|: u=\alpha(x \circ y)\beta\}\]
where infimum is taken over arbitrary decompositions with $x\in M_\infty(V)$,  $y\in M_\infty(W)$, $\alpha\in M_{n,\infty}$, and $\beta\in M_{\infty,n}$.
\end{prop}
\begin{pf}
Suppose that $u\in  M_n(V \otimes^s W)$ such that $\|u\|_s<1$. Let $\epsilon=1-\|u\|_s >0$. Let $u_k\in M_n(V\ot W)$ be the sequence such that $\|u-u_k\| \to 0$. We may assume that $\|u-u_1\|<\epsilon$ and $\|u_{k+1}-u_k\|\leq \frac{\epsilon}{2^{k+1}}$ for all $k\in \mathbb{N}$, from which we get $u=u_1+\ds\sum_{k=1}^{\infty}(u_{k+1}-u_k)$. Let $t_k= u_{k+1}-u_k$ and $t_0=u_1$. Then we have that $u=\ds\sum_{k=0}^{\infty}t_k$. As $t_k\in M_n(V\ot W)$, there exist $\alpha_k\in M_{n,p_k}$, $\beta\in M_{p_k,n}$, $x_k\in M_{p_k}(V)$, and $y_k\in M_{p_k}(W)$ such that
 $t_k=\alpha_k(x_k\circ y_k)\beta_k$ with $\|\alpha_k\|\|x_k\|\|y_k\|\|\beta_k\|< \|t_k\|+\frac{\epsilon}{2^{k+1}}$. We can assume that
 $\|x_k\|=\|y_k\|=1$, and $\|\alpha_k\|=\|\beta_k\|<(\|t_k\|+\frac{\epsilon}{2^{k+1}})^{\frac{1}{2}}$, for $k\in \mathbb{N}$, and
 $\|\alpha_0\|\|\beta_0\|< 1-\epsilon$. Then we have $\ds\sum_{k=0}^{\infty}\|\alpha_k\|\|\beta_k\|< 1$. Now choose an increasing sequence $c_k\geq 1$ with $c_k\to \infty $ such that $\ds\sum_{k=0}^{\infty}c_k\|\alpha_k\|\|\beta_k\|< 1$. Put $v=\bigoplus c_k^{-1}v_k\in M_{\infty}(V)$, $w=\bigoplus c_k^{-1}w_k\in M_{\infty}(V)$,  $\alpha=[c_1\alpha_1,c_2\alpha_2, \cdots,c_r\alpha_r\cdots ]\in M_{n,\infty}$, and $\beta=\alpha^{t}$. We then have $u=\alpha(v\circ w)\beta$ and all $v$, $w$, $\alpha$, and $\beta $ have norm $<1$.
\end{pf}

%
%
%

We  define the schur bounded map from operator space $V$ to operator space $W^*$ by using the identification $SB(V,W^*)=(V\otimes^s W)^*$, i.e. $SB(V,W^*)$ is isometrically isomorphic to $(V\otimes^s W)^*$.

\begin{prop}\label{sc2}
If $\varphi$ is a schur bounded map from $V$ to $W$. Then $\varphi^*: W^*\to V^*$ is also schur bounded with the same schur norm.
\end{prop}
\begin{pf}
Given that $\varphi$ is a schur bounded map from $V$ to $W$, meaning that $i\circ \varphi: V\to W^{**}$ is schur bounded, where $i:W\to W^{**}$ is the natural embedding. Thus  there exists $\phi\in  (V \otimes^s W^*)^* $ such that $\phi(v\otimes f)= i\circ \varphi(v)(f)=f(\varphi(v)) $, for all $v\in V$ and $f\in W^*$, with $\|\phi\|=\|i\circ \varphi\|$. Now consider the map $\theta : W^* \otimes^s V\to V \otimes^s W^*$, which is completely isometric by Proposition \ref{sc1}. So $\phi \circ \theta \in ( W^* \otimes^s V )^*$ with $\phi \circ \theta(f\otimes v)= \phi(v\otimes f)= f(\varphi(v))=\varphi^*(f)(v) $. Thus $\|\phi \circ \theta\|=\|\varphi^*\|$.   As $\theta $ is completely isometric map from $W^* \otimes^s V$ onto $V \otimes^s W^*$, so it follows that $\|\phi \circ \theta\|= \|\phi\|$. Thus $ \|i\circ \varphi\|= \|\varphi\|= \|\varphi^*\|$.
\end{pf}

If we use the following explicit definition of schur bounded map then it follows easily that every linear functional $f$ is schur bounded with $\|f\|_{sb}=\|f\|$.

\begin{defi}
Given operator spaces $V$ and $W$, let $\varphi: V  \to W^*$ be a linear map. We say that $\varphi$ is schur bounded  map if  $\|\varphi\|_{sb} = \displaystyle\sup\{\|[\varphi(x_{ij})(y_{ij})]\|: \|[x_{ij}]\|_{M_n(V)}\leq 1, \|[y_{ij}]\|_{M_n(W)}\leq 1, n\in \mathbb{N}\}< \infty$.
\end{defi}
Using the above definition, one can easily prove that $SB(V\times W,\C)=SB(V,W^*)$. Proposition \ref{sc2} can also be proved by using this explicit definition of schur bounded maps.


%
%
%
We now proceed to show that the schur tensor product of two matrix ordered space is matrix ordered. Recall that a complex vector space $V$ is said to be matrix ordered if (1) $V$ is a $^*$-vector space, (2) Each $M _n(V )$, $n \geq 1$, is partially ordered by a cone $M _n(V )_{+}\subseteq M _n(V )_{sa}$, the self adjoint part of $M_n(V)$, and (3) If $\gamma\in M_{m,n}$, then $\gamma^* M_m(V)^+\gamma\subseteq M_n(V)^+$. Also, by an involutive operator space, we mean   an operator space with an involution
such that, for each $n \in \mathbb{N}$ , $M_n(V)$ is an involutive Banach space with the natural involution,
i.e., the involution on $M_n(V)$ is an isometry.

\begin{prop}
Let $V$ and $W$ be an involutive operator spaces. If $u\in M_n(V \otimes_s W)_{sa}$ then $u$ has  a representation $u=\alpha(x\circ y)\alpha^*$,  where $\alpha\in M_{n,p}$, $x\in M_p(X)_{sa}$ and $y\in M_p(Y)_{sa}$. Moreover, $\|u\|_s=\inf\{\|\alpha\|^2\|\|x\|\|y\|: u=\alpha(x\circ y)\alpha^*, \alpha\in M_{n,p}, x\in M_p(X), y\in M_p(Y), p\in \mathbb{N}\}$.
\end{prop}
\begin{pf}
Let $u\in M_n(V \otimes_s W)_{sa}$ and $\epsilon>0$. Then there exist $\alpha\in M_{n,p}$, $\beta\in M_{n,p}$, $x\in M_p(V)$ and $y\in M_p(W)$ such that $u=\alpha(x\circ y)\beta $ with $\|u\|_{s}\leq \|\alpha\|\|x\|\|y\|\|\beta\|\leq \|u\|_{s}+\epsilon$. As $u$ is self adjoint, so we have  \\
\hspace*{3 cm} $u=\frac{1}{2}(u+u^*)$\\
\hspace*{3.34 cm} $=\frac{1}{2}(\alpha(x\circ y)\beta + \beta^*(x^*\circ y^*)\alpha^* )$.\\
\hspace*{3.44 cm}$=\left(
                  \begin{array}{cc}
                    \frac{\lambda \beta^*}{\sqrt{2}} & \frac{\lambda^{-1}\alpha}{\sqrt{2}} \\
                  \end{array}
                \right)
\left(
                                                      \begin{array}{cc}
                                                        0 & x^* \\
                                                        x & 0\\
                                                      \end{array}
                                                    \right)
                                                    \circ \left(
                                                             \begin{array}{cc}
                                                               0& y^* \\
                                                               y & 0\\
                                                             \end{array}
                                                           \right)\left(
                   \begin{array}{c}
                     \frac{\lambda  \beta}{\sqrt{2}} \\
                     \frac{\lambda^{-1}\alpha^*}{\sqrt{2}}\\
                   \end{array}
                 \right)
$ for any $\lambda> 0$.  Let $v:=\left(
                                                      \begin{array}{cc}
                                                        0 & x^* \\
                                                        x & 0\\
                                                      \end{array}
                                                    \right)$, $w:= \left(
                                                             \begin{array}{cc}
                                                               0& y^* \\
                                                               y & 0\\
                                                             \end{array}
                                                           \right)$, and $\tilde{\alpha}= \left(
                  \begin{array}{cc}
                    \frac{\lambda \beta^*}{\sqrt{2}} & \frac{\lambda^{-1}\alpha}{\sqrt{2}} \\
                  \end{array}
                \right)$. Then we have $u= \tilde{\alpha}(v\circ w)\tilde{\alpha}^*$ and $\|u\|_{s}\leq  \|\tilde{\alpha}\|^2\|v\|\|w\| \leq[\frac{1}{2}(\lambda^2\|\beta\|^2+\lambda^{-2}\|\alpha\|^2)]\|v\|\|w\|$, where $v$ and $w$ are self adjoint elements.
Now, by using the fact that $\ds\min_{\lambda>0}\frac{1}{2}(\lambda^2\|\beta\|^2+\lambda^{-2}\|\alpha\|^2)= \|\beta\|\|\alpha\|$, given $\delta>0$ choose $\lambda_0>0$ such that $\|\beta\|\|\alpha\|+\delta > \frac{1}{2}(\lambda_0^2\|\beta\|^2+\lambda_0^{-2}\|\alpha\|^2)$. We then have  $\|u\|_s\leq \|\tilde{\alpha}\|^2\|v\|\|w\|\leq \|\beta\|\|\alpha\|\|x\|\|y\|$. Thus, we get the desired norm condition.
\end{pf}

Based on the above,  we define  $(V \otimes^s W)_{+}$=cl$_{s}\{ \alpha (v \circ w)\alpha^*: v\in M_p(V)^{+},  w\in M_p(W)^{+}, \alpha\in M_{n,p} \}$, and we have the following, which can be proved easily.

\begin{prop}
For matrix  ordered operator spaces $V$ and $W$, $V \otimes^s W$ is a matrix ordered operator space.
\end{prop}
Let $X$ and $Y$ be operator spaces. For $f \in CB(X,M_p)$ and $g \in CB(Y,M_p)$, define $f\bar{\circ} g: X\ot Y\to M_p$ on elementary tensor as $f\bar{\circ} g(x\ot y)=f(x)\circ g(y)$ for $x\in X$ and $y\in Y$.  Define a norm on $M_n(X\ot Y)$ as $\|u\|_{s'} = \sup\{\|(f\bar{\circ} g)_n(u)\|\}$, where the supremum is taken over all $f \in M_p(X^*)_1$
and $g \in M_p( Y^*)_1$, $p\in \mathbb{N}$. Let $X\ot^{s'} Y$ denote the completion of $X\ot Y$ in the  $\|\cdot\|_{s'}$-tensor norm.

\begin{prop}\label{sc22}
For operator spaces $X$ and $Y$,  the  natural embedding $\theta: X\ot^{s'} Y\to SB(X^*\times Y^*,\mathbb{\C})$ is isometric.
\end{prop}
\begin{pf}
The canonical map $\theta$ is determined by $\theta(x\ot y)(f,g)=f(x)g(y)$ for $x\in X$, $y\in Y$, $f\in X^*$ and $g\in Y^*$. For any $u=\ds\sum_{t=1}^{k}a_t\ot b_t\in X\ot Y$, by definition, $\|\theta(u)\|_{sb}=\sup\{\|[\theta(u)(f_{ij},g_{ij})]\|: [f_{ij}]\in M_n(X^*)_1, [g_{ij}]\in M_n(Y^*)_1, n\in \mathbb{N}\}= \sup\{\|[\ds\sum_{t=1}^k f_{ij}(a_t)g_{ij}(b_t)]\|: [f_{ij}]\in M_n(X^*)_1, [g_{ij}]\in M_n(Y^*)_1, n\in \mathbb{N}\}=\sup\{\|[ [f_{ij}]\bar{\circ} [g_{ij}](u)]\|: [f_{ij}]\in M_n(X^*)_1, [g_{ij}]\in M_n(Y^*)_1, n\in \mathbb{N}\}$, by using the identification $M_n(X^*)=CB(X,M_n)$, which is same as $\|u\|_{s'}$.
\end{pf}

Next we consider the useful variation of the last proposition, which shows that the dual of the schur tensor norm is  the $\|\cdot\|_{s'}$-norm.

\begin{prop}\label{sc3}
For operator spaces $X$ and $Y$, the   natural embedding $\psi: X^*\ot^{s'}Y^*\to SB(X\times Y,\mathbb{\C})$ is isometric.
\end{prop}
\begin{pf}
We have to show that $\|\ds\sum_{t=1}^k f_t\ot g_t\|_{s'}=\sup\{\|[\ds\sum_{t=1}^k f_t(x_{ij})g_t(y_{ij})]\|: [x_{ij}]\in M_n(X)_1, [y_{ij}]\in M_n(Y)_1, n\in \mathbb{N}\}$. By  Proposition \ref{sc22}, we have an isometric map  $X^*\ot^{s'} Y^*\to SB(X^{**}\times Y^{**},\mathbb{\C})$. Therefore, the right-hand side is dominated by $\|u\|_{s'}$ for $u=\ds\sum_{t=1}^k f_t\ot g_t$. Now let $ [F_{ij}]\in M_n(X^{**})_1$, $[G_{ij}]\in M_n(Y^{**})_1$ then there exist $\widehat{x_\lambda}\in M_n(X) $ and $\widehat{y_\nu}\in M_n(Y)$ such that $\widehat{x_\lambda}$ converges to $[F_{ij}]$ and
 $\widehat{y_\nu}$ converges to $[G_{ij}]$ in the point-norm topology by (~\cite{effros}, Proposition 4.2.5). Thus equality holds.
\end{pf}

The identification  $M_n(SB(X\times Y,\mathbb{\C}))=SB(X\times Y,M_n)$ endows $SB(X\times Y,\mathbb{C})$ with an operator space structure and make  the map, defined in Proposition \ref{sc3}, completely isometric. Also note that, like the projective and injective norm, schur and delta norm are in perfect duality  in the finite dimensional setting.

\section{Equivalence of the Schur norm}
Let $E\subseteq A$ be an operator subspace of a $C^*$-algebra $A$. Then $E$ is said to be completely complemented if there is a completely bounded (cb) projection $P$ from $A$ onto $E$. In analogy to the operator space projective tensor product, schur tensor does not respect subspaces in general but  behaves well for completely complemented subspaces:

\begin{lem}\label{sc6}
Let $E$, $F$ be completely  complemented subspaces of the $C^*$-algebras $A$ and $B$ complemented by cb projection having  cb norm 1, respectively. Then $E\otimes^s F$ is a closed subspace of $A\otimes^s B$.
\end{lem}
\begin{pf}
By an assumption, there are cb projections  $P$ from  $A$ onto  $E$, and $Q$ from $ B$ onto $F$ with $\|P\|_{cb}=1$, and $\|Q\|_{cb}=1$.  Therefore, by Corollary \ref{sc113}, $P\ot Q: A\otimes^s B \to E\otimes^s F$ is a bounded map and $\|P\ot Q\|\leq 1$. Now, for $u\in  E\otimes F$, $P\ot Q(u)=u$, giving that $\|u\|_{E\otimes^s F}\leq \|u\|_{A\otimes^s B}$. Hence $E\otimes^s F$ is a closed subspace of $A\otimes^s B$.
\end{pf}

In particular, if  $E$ and $F$ are finite dimensional $C^*$-subalgebras of the $C^*$-algebras $A$ and $B$, respectively. Then $E\otimes^s F$ is a closed $^*$-subalgebra of $A\otimes^s B$ by (~\cite{bal}, II 6.10.4(iii)). Also, for von Neumann algebras $M$ and $N$, $Z(M)\ot^s Z(N)$ is a closed $^*$-subalgebra of $M\ot^s N$ by (~\cite{kadison}, \S 3, Theorem C).

Recall that the tracially bounded norm on the algebraic tensor product of two $C^*$-algebras $A$ and $B$ is defined as
\[\|u\|_{tb}=\inf\{\ds\sum_{k=1}^N \|[a_{ij}^k]\|\|[b_{ji}^k]\|:u=\ds\sum_{k=1}^N  n^{-1}\ds\sum_{i,j=1}^n a_{ij}^k\ot b_{ji}^k \}.\]
\begin{lem}
 $\|\cdot\|_s\leq \|\cdot\|_{tb}$ on $A\otimes B$.
\end{lem}
\begin{pf}
Let $u=n^{-1}\ds\sum_{k=1}^{N} \ds\sum_{i,j=1}^n a_{ij}^{k}\ot b_{ji}^{k}$, then $\|u\|_s \leq n^{-1} \ds\sum_{k=1}^{N} \|\ds\sum_{i,j=1}^n a_{ij}^{k}\ot b_{ji}^{k}\|_s=
n^{-1}\ds\sum_{k=1}^{N} \|\left(\begin{matrix}
               1 \; 1 \;\;\ldots 1
  \end{matrix}\right) [a_{ij}^{k}]\circ [b_{ji}^{k}]\left(\begin{matrix}
  1  \;1 \;\; \ldots 1
  \end{matrix}\right)^{t} \|_s$
$\leq n^{-1} \ds\sum_{k=1}^{N} \|\left(\begin{matrix}
               1 \; 1 \;\;\ldots 1
  \end{matrix}\right)\|$ $\|[a_{ij}^{k}]\|\| [b_{ji}^{k}]\|\|\left(\begin{matrix}
  1  \;1 \;\; \ldots 1
  \end{matrix}\right)^{t} \|= \ds\sum_{k=1}^{N} \|[a_{ij}^{k}]\| \| [b_{ji}^{k}]\| $, and so $\|u\|_s\leq \|u\|_{tb}$.
\end{pf}

Therefore, we have the following comparison between the various tensor norms:\\
$\|\cdot\|_\lambda \leq \|\cdot\|_{s'} \leq\|\cdot\|_{\min} \leq \|\cdot\|_{\max} \leq \|\cdot\|_h \leq \|\cdot\|_{\wedge} \leq \|\cdot\|_{s} \leq \|\cdot\|_{tb} \leq \|\cdot\|_{\gamma} $.\\

We now look at the equivalence of the schur tensor norm with these norms.

\begin{lem}\label{sc7}
If $n\in \mathbb{N}$ then in $M_n\ot M_n$,
\[\|\ds\sum_{j=1}^n e_{j1}\ot e_{j1}\|_s= n^{1/2}.\]
\end{lem}
\begin{pf}
Let $\al$ be a $1 \times n$ matrix with 1 in the (1,1) position and all other entries are zeros, and $\beta$ be a $n\times 1$ matrix with 1 in all entries. Let $x$ and $y$ be $n\times n$ matrices in $M_n(M_n)$ with first row $e_{11}, e_{21},...,e_{n1}$, and all other entries are zeros. Now it follows, from  $C^*$-identity and Ruan's axioms of operator space,  that $\|\al\|=1$, $\|\beta\|=n^{1/2}$, and $\|x\|=\|y\|=1$. Since $\ds\sum_{j=1}^n e_{j1}\ot e_{j1}= \al(x \circ y)\beta$, so $\|\ds\sum_{j=1}^n e_{j1}\ot e_{j1}\|_s\leq n^{1/2}$. Other inequality is obvious  by (~\cite{Ajay}, Lemma 3.3(i)) and the fact that $\|\cdot \|_{\wedge}\leq \|\cdot \|_{s}$.
\end{pf}

\begin{lem}\label{sc12}
Let $M$ and $N$ be von Neumann algebras and $T : M\times N \to \C$ be
a separately normal bilinear form. Then, for each $n \in N$, the map $T_n : M_n(M) \times M_n(N) \to M_n$ given by
$T_n((a_{ij}), (b_{ij})) = (T(a_{ij} , b_{ij}))$ is separately normal.
\end{lem}
\begin{pf}
We  only show that $T_n$, $n\in \mathbb{N}$ fixed, is normal in the first variable. In the second variable, result follows on the similar lines. Let $(a_{\lambda})$ be an increasing net of positive elements in  $M_n(M)$ such that $a_{\lambda}$ is $w^*$-convergent to $a\in M_n(M)$. Let $b=[b_{ij}]$ be fixed matrix in $ M_n(N)$. Since $a_{\lambda}\in M_n(M)$, so let $a_{\lambda}= [a_{ij}^{\lambda}]$ and $a=[a_{ij}]$. Since $T$ is separately normal, so $T(a_{ij}^{\lambda}, b_{ij})$ is $w^*$-convergent to $T(a_{ij}, b_{ij})$ for each $i,j$. Thus $(T(a_{ij}^{\lambda}, b_{ij}))$ is $w^*$-convergent to $(T(a_{ij}, b_{ij}))$, showing that $T_{n}(a_{\lambda}, b)$ is $w^*$-convergent to  $T_n(a,b)$.
\end{pf}
\begin{prop}\label{sc9}
Let $A$ and $B$ be $C^*$-algebras and $\p: A \times B \to \C$ be a schur bounded bilinear form. Then there exists a unique separately normal schur  bounded bilinear form $\tilde{\p}: A^{**}\times B^{**}\to \C$ such that $\|\p\|_{sb}=\|\tilde{\p}\|_{sb}$.
\end{prop}
\begin{pf}
Since $\p: A \times B \to \C$ is a schur bounded bilinear form, it is in particular bounded bilinear form and thus determines a unique separately normal bilinear form  $\tilde{\p}: A^{**} \times B^{**} \to \C$ with $\|\p\|=\|\tilde{\p}\|$ by (~\cite{haag}, Corollary 2.4). We show that $\|\p\|_{sb}=\|\tilde{\p}\|_{sb}$. For  $n\in \mathbb{N}$, consider the map $\tilde{\p}_n: M_n(A^{**}) \times M_n(B^{**}) \to M_n$ defined as $\tilde{\p}_n((a_{ij}), (b_{ij})) = (\tilde{\p}(a_{ij} , b_{ij}) )$. Let $a^{**}\in M_n(A^{**})$ and $b^{**}\in M_n(B^{**})$ with
$\|a^{**}\|\leq 1$ and $\|b^{**}\|\leq 1$. Since the unit ball of  $M_n(A)$ is $w^*$-dense in the unit ball of  $M_n(A^{**})$, so we obtain  a net $(a_{\lambda})$ (resp., $(b_\nu)$) in $M_n(A)$ (resp., $M_n(B)$) which
is $w^*$-convergent to $a^{**}$ (resp., $b^{**}$) with $\|a_{\lambda}\|\leq 1$ (resp., $\|b_{\nu}\|\leq 1$). By Lemma \ref{sc12},  $\tilde{\p}_n$  is separately normal so  $\|\tilde{\p}_n(a^{**}, b^{**})\|\leq \lim \inf_{\lambda, \nu}\|\tilde{\p}_n(\widehat{a_{\lambda}},\widehat{b_{\nu}} )\|= \lim \inf_{\lambda, \nu}\|\p_n(a_{\lambda},b_{\nu} )\|\leq \|\p_n\| $. Thus $\|\tilde{\p}_n\|\leq \|\p_n\|\leq \|\p\|_{sb}$  for every $n\in \mathbb{N}$.  Clearly, $\|\p\|_{sb}\leq \|\tilde{\p}\|_{sb}$ as $\p$ being the restriction of $\tilde{\p}$. Hence $\|\p\|_{sb}=\|\tilde{\p}\|_{sb}$.

\end{pf}
\begin{cor}
For $C^*$-algebras $A$ and $B$, the embedding $\nu= i_{A}\ot^{s} i_{B}$ of $A\otimes^s B$ into $A^{**}\otimes^s B^{**}$ is an isometry.
\end{cor}

For $C^*$-algebras $A$ and $B$, we know that  $SB(A,B^*)=(A\otimes^s B)^*$. So, by Proposition \ref{sc9}, we have an isometry $\alpha: (A\otimes^s B)^* \to  (A^{**}\otimes^s B^{**})^*$ given by
$(\alpha\phi)(a^{**}\otimes b^{**})= \Phi(a^{**})(b^{**})$.  Let $\omega:= \alpha^{*} \circ i: A^{**}\otimes^s B^{**} \to (A\otimes^s B)^{**} $, where $i$ is the natural embedding from $A^{**}\otimes^s B^{**}$ into  $(A^{**}\otimes^s B^{**})^{**}$. By the definition, it is clear that $\omega$ is norm reducing.

\begin{prop}\label{sc14}
For $C^*$-algebras $A$ and $B$, we have the following :\\
\emph{(1)} $\|\cdot\|_s \approx \|\cdot\|_\gamma $ if and only if either $A$ or $B$ is subhomogeneous.\\
\emph{(2)} $\|\cdot\|_s \approx \|\cdot\|_{tb} $ if and only if either $A$ or $B$ is subhomogeneous.\\
\emph{(3)} $\|\cdot\|_\wedge \approx \|\cdot\|_{tb} $ if and only if either $A$ or $B$ is subhomogeneous.\\
\emph{(4)} $\|\cdot\|_s \approx \|\cdot\|_h $ if and only if $A$ or $B$ is finite dimensional or $A$ and $B$ are infinite dimensional subhomogeneous.\\
\emph{(5)} $\|\cdot\|_s \approx \|\cdot\|_{\max} $ if and only if $A$ or $B$ is finite dimensional.\\
\emph{(6)} $\|\cdot\|_s \approx \|\cdot\|_{s'} $ if and only if $A$ or $B$ is finite dimensional. \\
 \emph{(7)} If either $A$ or $B$ is subhomogeneous then $\|\cdot\|_s \approx \|\cdot\|_{\wedge} $. \\
\end{prop}

\begin{pf}
(1): Suppose that  $A$ is subhomogeneous. Then $\|\cdot\|_\wedge\approx \|\cdot\|_\gamma $ by (~\cite{Ajay}, Theorem 7.2), which in turn implies that $\|\cdot\|_s \approx \|\cdot\|_\gamma $.

Conversely, suppose that $\|u\|_\gamma \leq K\|u\|_s$,  for some constant $K> 0$, for all $u\in A\otimes B$. Since the map $\omega$ is norm reducing, so $\|\omega(u)\|_{**}\leq \|u\|_s$ for all $u\in A^{**}\otimes B^{**}$. Also, by  (~\cite{vandee3}, Theorem 2.3), we have  $\|u\|_\gamma \leq 2\|\mu (u)\|_{**}$ for all $u\in A^{**}\otimes B^{**}$, where $\mu$ is a map from $A^{**}\otimes^\gamma B^{**}$ into $(A\otimes^\gamma B)^{**}$, and $\|\cdot\|_{**}$ denotes the relevant second dual norm. Therefore, $\|u\|_\gamma \leq 2\|\mu (u)\|_{**} \leq 2K \|\omega (u)\|_{**}\leq 2K \|u\|_s$
for all $u\in A^{**}\otimes B^{**}$. Now suppose that neither $A$ nor $B$ are subhomogeneous. Then, for some positive integer $n$, $A^{**}$  and $B^{**} $ contain a copy of $M_n$, which by Lemma \ref{sc6} implies that $M_n \otimes^s  M_n$ embeds isometrically into  $A^{**}\otimes^s B^{**}$. Also, $M_n \otimes_\gamma  M_n$ embeds isometrically into  $A^{**}\otimes_\gamma  B^{**}$ as there is a conditional expectation  from $A^{**}$ onto $M_n$. By Lemma \ref{sc7} and (~\cite{Ajay}, Lemma 3.3), it follows that
\[n= \|\ds\sum_{j=1}^n e_{j1}\ot e_{j1}\|_\gamma \leq 2K \|\ds\sum_{j=1}^n e_{j1}\ot e_{j1}\|_s= 2K n^{1/2}.\] Thus one of  $A^{**}$ and  $B^{**}$  cannot contain a type $I_n$ factor for $n \geq 4 K^{2}$,  giving that either $A$ or $B$ is subhomogeneous for $n\geq 4 K^{2}$.\\
(2): Suppose that $\|\cdot\|_s \approx \|\cdot\|_{tb} $.  So, by (~\cite{trbl}, Theorem 1),  $\|\cdot\|_s \approx \|\cdot\|_\gamma $. Thus either $A$ or $B$ is subhomogeneous by (1).

Conversely, suppose that either $A$ or $B$ is subhomogeneous. Then $\|\cdot\|_s \approx \|\cdot\|_\gamma $ by (1), and so  $\|\cdot\|_s \approx \|\cdot\|_{tb} $ by (~\cite{trbl}, Theorem 1).\\
(3): follows as in (2).\\
(4): follows directly from (~\cite{Ajay}, Theorem 6.1 and Theorem 7.4).\\
(5): Since $\|\cdot\|_s \approx \|\cdot\|_{\max} $ implies that $\|\cdot\|_h \approx \|\cdot\|_{\max} $. Hence the result follows from (~\cite{Itoh}, Theorem). Converse is obvious.\\
(6): Suppose that $A$ and $B$ are infinite dimensional. Choose  maximal abelian subalgebras $A_1$ and $B_1$ inside $A$ and $B$, respectively, by Zorn's Lemma. By (~\cite{ringrose}, Exercise 4.6.12), $A_1$ and $B_1$ are infinite dimensional.  Since $\|\cdot\|_s \approx \|\cdot\|_{s'} $, so $\|\cdot\|_h\approx \|\cdot\|_{s'}$. Since $\|\cdot\|_{h}$ and $\|\cdot\|_{s'}$ are injective, so $\|\cdot\|_h\approx \|\cdot\|_{s'}$ on $A_1\ot B_1$. As $A_1$ and $B_1$ are commutative, so   $\|\cdot\|_h\approx \|\cdot\|_{\max}$ on $A_1\ot B_1$. Thus $A_1$ or $B_1$ is finite dimensional by (~\cite{Itoh}, Theorem), a contradiction. Converse is trivial.\\
(7) follows from (3). Note that if $A$ or $B$ is finite dimensional then $\|\cdot\|_s \approx \|\cdot\|_h $ by (4) which in turn implies $\|\cdot\|_s \approx \|\cdot\|_\wedge $. But in this case by (~\cite{Ajay}, Theorem 6.1), we have that $\|\cdot\|_s \approx \|\cdot\|_\gamma $. This implies that $A$ or $B$ is subhomogeneous, which is not true in general.
\end{pf}

Proposition \ref{sc14} leaves open the interesting question whether $\|\cdot\|_s \approx \|\cdot\|_{\wedge} $ implies  $A$ or $B$ is subhomogeneous or not.


%
%

\begin{prop}
For operator spaces $V$ and $W$, $\max(V)\widehat{\ot}W \cong max(V)\ot^{s}W$.
\end{prop}
\begin{pf}
Since $JCB(max(V)\times W,\C)\subseteq SB(max(V)\times W,\C)\subseteq B(max(V)\times W,\C)$, so (~\cite{effros}, \S 3.3.9) implies that  $JCB(max(V)\times W,\C)= SB(max(V)\times W,\C)$, and hence  the result.
\end{pf}

We now summarize the  norms  of various elements in $M_n\ot M_n$, the calculation of which can be carried out as in (~\cite{Ajay}, Lemma 3.1) and Lemma \ref{sc7}. The details are left to the reader.\\

\begin{tabular}{|r|r|r|r|r|r|r|}
\hline
\text{Elements}& $\|\cdot\|_{\min}$ & $\|\cdot\|_h$ & $\|\cdot\|_{\wedge}$ & $\|\cdot\|_{s}$ & $\|\cdot\|_{tb}$ & $\|\cdot\|_{\gamma}$ \\
\hline
$\ds\sum_{j=1}^n e_{1j}\ot e_{jj}$&1& $\sqrt{n}$&$\sqrt{n}$&$\sqrt{n}$&$\sqrt{n}$&$\sqrt{n}$\\
\hline
$\ds\sum_{j=1}^n e_{j1}\ot e_{jj}$ &1 & 1&$\sqrt{n}$&$\sqrt{n}$&$\sqrt{n}$&$\sqrt{n}$\\
\hline
$\ds\sum_{j=1}^n e_{1j}\ot e_{1j}$ &$\sqrt{n}$& $\sqrt{n}$&$\sqrt{n}$&$\sqrt{n}$&$n$&$n$\\
\hline
$\ds\sum_{j=1}^n e_{j1}\ot e_{j1}$ &$\sqrt{n}$& $\sqrt{n}$&$\sqrt{n}$&$\sqrt{n}$&$n$&$n$\\
\hline
$\ds\sum_{i,j=1}^n e_{ij}\ot e_{ij}$ &$n$& $n$&$n$&$n$&$n$&$n$\\
\hline
$\ds\sum_{i,j=1}^n e_{i1}\ot e_{ij}$ &$n$& $n$&$n$&$n$&$n^{\frac{3}{2}}$&$n^{\frac{3}{2}}$\\
\hline
$\ds\sum_{i,j=1}^n e_{ii}\ot e_{ij}$ &$\sqrt{n}$& $n$&$n$&$n$&$n$&$n$\\
\hline
$\ds\sum_{j=1}^n e_{ji}\ot e_{ij}$ &1& $n$&$n$&$n$&$n$&$n$\\
\hline
$\ds\sum_{j=1}^n e_{1j}\ot e_{j1}$ &1& $n$&$n$&$n$&$n$&$n$\\
\hline
$\ds\sum_{i,j=1}^n e_{jj}\ot e_{ij}$ &$\sqrt{n}$& $\sqrt{n}$&$n$& $n$&$n$&$n$\\
\hline
$\ds\sum_{j=1}^n e_{1j}\ot e_{ij}$&$n$& $n$&$n$&$n$&$n^{\frac{3}{2}}$&$n^{\frac{3}{2}}$\\
\hline
\end{tabular}
\newpage
\begin{tabular}{|r|r|r|r|r|r|r|}
\hline
\text{Elements}& $\|\cdot\|_{\min}$ & $\|\cdot\|_h$ & $\|\cdot\|_{\wedge}$ & $\|\cdot\|_{s}$ & $\|\cdot\|_{tb}$ & $\|\cdot\|_{\gamma}$ \\
\hline
$\ds\sum_{j=1}^n e_{1j}\ot e_{ji}$&$\sqrt{n}$& $n^{\frac{3}{2}}$&$n^{\frac{3}{2}}$&$n^{\frac{3}{2}}$&$n^{\frac{3}{2}}$&$n^{\frac{3}{2}}$\\
\hline
$\ds\sum_{j=1}^n e_{j}\ot e_{jj}\in l_n^\infty\ot M_n$ &1 & $\sqrt{n}$&$\sqrt{n}$&$\sqrt{n}$&$\sqrt{n}$&$\sqrt{n}$\\
\hline
\end{tabular}


\begin{thebibliography}{a}
\bibitem{bal}Blackadar, B. , Operator algebras: Theory of $C^{*}$-algebras and von Neumann algebras, \textit{Springer-Verlag Berlin Heidelberg}, 2006.
    \bibitem{geom}  Blecher, D. P., Geometry of the tensor product of $C^*$-algebras,  \textit{Math. Proc. Camb. Phil. Soc.} 104 (1988), 119--127.
    \bibitem{trbl} Blecher, D. P., Tracially completely bounded multilinear maps on $C^*$-algebras, \textit{J. London Math. Soc.} s2-39 (1989), 514--524.
\bibitem{blecher}  Blecher, D. P. and  Paulsen,V. I., Tensor Products of operator Spaces, \textit{J. Func. Anal. }99 (1991), 262--292.


\bibitem{ER} Effros, E.G. and Ruan, Z.-J., A new approach to opertor spaces, \textit{Canad. Math. Bull. }34 (1991), 329-337.
\bibitem{effros} Effros, E. G. and Ruan, Z. J., Operator spaces, \textit{Claredon Press-Oxford}, 2000.
\bibitem{kadison} Ge., L. and Kadison, R., On tensor product of von Neumann algebras, \textit{ Invent. math} 123 (1996),453--466.
\bibitem{haag} Haagerup, U., The Grothendieck inequality for bilinear forms on $C^{*}$-algebras, \textit{ Adv. Math. 56} (1985), 93--116.
\bibitem{Itoh1} Itoh, T., On the completely bounded maps of a $C^*$ algebra to
    its dual space,  \textit{ Bull. London Math. Soc. } (1987) 19, 546--550.
\bibitem{Itoh} Itoh, T., The maximal $C^*$-norm and the Haagerup norm, \textit{Math. Proc. Camb. Phil. Soc.}  107 (1990), 109--114.
\bibitem{haageritoh} Haagerup, U. and  Itoh, T.,   Grothendieck type norms for bilinear forms on $C^*$-algebras,
\textit{ J. Operator Theory}  34  (1995), 263--283.
\bibitem{itoh} Itoh, T, Completely positive decompositions from duals of $C^*$-algebras to von Neumann algebras,
\textit{ Math. Japonica} 51 (2000), 89--98.
    \bibitem{r2}Jain, R. and Kumar, A., Operator space projective tensor product: Embedding into second dual and ideal structure, \textit{ To appear in Proc. Edin. Math. Soc}, Available on arXiv:1106.2644v1 [math.OA].
 \bibitem{ringrose} Kadison, R. V. and  Ringrose, J. R., Fundamentals of the theory of operator algebras I. \textit{Academic Press}, 1983.
\bibitem{Ajay} Kumar, A. and  Sinclair, A. M., Equivalence of norms on operator space tensor products of $C^{*}$-algebras, \textit{Trans. Amer. Math. Soc.} 350  (1998), 2033--2048.
\bibitem{vandee3}Kumar, A. and Rajpal, V., Projective tensor product of $C^{*}$-algebras, Available on  arxiv:1305.0791v1 [math.OA].
\bibitem{laur}Laursen, Kjeld B., Tensor products of Banach algebras with involution, \textit{Trans. Amer. Math. Soc.} 136 (1969), 467--487.
\end{thebibliography}
\end{document}